\documentclass{article}
\usepackage{comment}
\usepackage{graphicx}
\usepackage{indentfirst}
\usepackage{url}
\usepackage{color}
\usepackage{amsmath}
\usepackage{amsfonts}
\usepackage{amsthm}
\usepackage{amssymb}
\usepackage{bm}
\usepackage{enumitem}
\begin{document}

\newtheorem{theorem}{Theorem}[section]
\newtheorem{lemma}[theorem]{Lemma}

\theoremstyle{definition}
\newtheorem{definition}[theorem]{Definition}
\newtheorem{example}[theorem]{Example}
\newtheorem{xca}[theorem]{Exercise}

\theoremstyle{remark}
\newtheorem{remark}[theorem]{Remark}
\theoremstyle{conjecture}
\newtheorem{conjecture}[theorem]{Conjecture}

\abovedisplayskip 6pt plus 2pt minus 2pt \belowdisplayskip 6pt
plus 2pt minus 2pt
%%%%%%%%%%%%%%%%
\def\vsp{\vspace{1mm}}
\def\th#1{\vspace{1mm}\noindent{\bf #1}\quad}
\def\proof{\vspace{1mm}\noindent{\it Proof}\quad}
\def\no{\nonumber}
\newenvironment{prof}[1][Proof]{\noindent\textit{#1}\quad }
{\hfill $\Box$\vspace{0.7mm}}
\def\q{\quad} \def\qq{\qquad}
%\allowdisplaybreaks[4]
%%%%%%%%%%%%%%%%%%%%%%%%%%%%%%%%%%%%%%%%%%%%%%%%%%%%%%%%%%%%%%%%%%%%%%%%%%%%%%%%%%%%%%%%%%%%%%%
%%-------------------       Beginning of  Author's Definitions       -------------------%%
%%                     Note: You may add your own definitions here.
\def\IDiscrim{{\rm IDiscrim}}
\def \Res  {{\rm Res}}
\def \MCproj  {{\rm MCproj}}
\def \cs  {{\rm cs}}
\def \cso  {{\rm cso}}
\def \ms  {{\rm ms}}
\def \mso  {{\rm mso}}
\def \discrim  {{\rm Disc}}
\def \sqrfree  {{\rm sqrfree}}
\def \sqk  {{\rm sqk}}
\def \lc  {{\rm lc}}
\def \Bproj  {{\tt IR}}
\def \Bp  {{\tt IR}}
\def  \zero {{\rm Zero}}
\def  \Nproj {{\tt Np}}
\def  \NTproj {{\tt NTp}}
 \def  \AlgNproj {{\tt Nproj}}
\def  \Hproj {{\tt Hp}}
\def  \Hp {{\tt Hp}}
\def  \Npl {{\tt NpL}}
\def  \Nfproj {{\tt Nfp}}
\def \Proineq {{\tt DPS}}
\def \ProCop {{\tt ProCop}}
\def \Bprojection {{\tt Bproj}}
\def \HNprojection {{\tt HNproj}}
\def  \HNproj {{\tt HNp}}
\def \Nprojection {{\tt Nproj}}
\def \Hprojection {{\tt Hproj}}
\def \Findk {{\tt Findk}}
\def \SRes {{\tt SRes}}
\def \Findinf {{\tt Findinf}}
\def \m  {{\rm m}}
\def \RR {{\mathbb R}}
\def \QQ {{\mathbb Q}}
\def \CC {{\mathbb C}}
\def  \Zar {{\tt Zar}}
\def  \Card {{\tt Card}}
\def \ZZ {{\mathbb Z}}
\def \AA {{\mathbb A}}
\newcommand{\va}{\bm{\alpha}}
\newcommand{\vb}{\bm{\beta}}
\newcommand{\xx}{\bm{x}}
\newcommand{\yy}{\bm{y}}
\newcommand{\zz}{\bm{z}}
\newcommand{\XX}{\bm{X}}
\newcommand{\YY}{\bm{Y}}

\def \R {{\mathbb R}}
\def \C {{\mathbb C}}
\def \rr {{\mathcal R}}
\def \RAGlib {{\tt RAGlib}}
\def \FI {{\tt FI}}
\def \QEPCAD {{\tt QEPCAD}}
\def \PCAD {{\tt PCAD}}
\def \TwoPro {{\tt PSD-HpTwo}}
\def \TwoHp {{\tt HpTwo}}
\def \RAGMaple {{\tt RAGMaple}}
\def \Raglib {{\tt Raglib}}
\def \FI {{\tt FI}}
\def \QEPCAD {{\tt QEPCAD}}
\def \PCAD {{\tt PCAD}}
\def \TCPT {{\tt CMT}}
\def \TCPTC {{\tt TCPTC}}
\def \OO {{\mathcal{O}}}

%%-------------------         the end of  Author's Definitions           -------------------%%

%\AuthorMark{Jingjun Han}                             %%%  appear on the head of even pages  %%%

%\TitleMark{Multivariate Discrim. and Iterated Resultant}  %%% Running Title, appear on the head of odd pages  %%%

\title{Multivariate Discriminant and Iterated Resultant       %%%   Main Title of your paper  %%%
\footnote{Supported by China Scholarship Council}}                  %%%   the Fund which you are supported by  %%%

\author{Jingjun Han\footnote{School of Mathematical Sciences $\&$ Beijing International Center for Mathematical Research, Peking University, Beijing 100871, China, hanjingjun@pku.edu.cn}}
\maketitle

\textbf{Abstract}  In this paper, we study the relationship between iterated resultant and multivariate discriminant. We show that, for generic form $f(\xx_n)$ with even degree $d$, if the polynomial is squarefreed after each iteration, the multivariate discriminant $\Delta(f)$ is a factor of the squarefreed iterated resultant. In fact, we find a factor $\Hp(f,[x_1,\ldots,x_n])$ of the squarefreed iterated resultant, and prove that the multivariate discriminant $\Delta(f)$ is a factor of $\Hp(f,[x_1,\ldots,x_n])$. Moreover, we conjecture that
$\Hp(f,[x_1,\ldots,x_n])=\Delta(f)$ holds for generic form $f$, and show that it is true for generic trivariate form $f(x,y,z)$.      % the abstract

%Cylindrical Algebraic Decomposition, semi-definiteness, polynomial, resultant, multivariate discriminant.% the keywords

%\MRSubClass{13P15,  68W30}      % MR(2000) Subject Classification

%\baselineskip 15pt

\section{Introduction}% with even degree
In this paper, for generic form $f(\xx_n)$, we study the relationship between multivariate discriminant and iterated resultant.
This topic has a long story. In 1868, \cite{henrici1866certain} considered the problem of decomposing two times iterated discriminant (the discriminant of the discriminant) of a generic trivariate form $f(x,y,z)$, and found it has a natural
factorization of the shape $\discrim(\discrim(f(1,y,z),z),y)=cPQ^2R^3$. In particular, the polynomial $P$ is the multivariate discriminant of the generic trivariate form. %he gave the expected factorization of such a repeated discriminant.
Recently, \cite{buse2009explicit} showed that each of the factors are irreducible. For related works, see for example, \cite{mccallum1999factors,lazard2009iterated}.

However, all of these previous results are about two times iterated discriminant. Our goal is to extend these results to cases of arbitrary iterated times, that is, for more than two times iteration. One of the main results of this paper is Theorem \ref{thm:main}, it claims that for generic form $f(\xx_n)$ with even degree $d$, if the polynomial is squarefreed after each iteration, the multivariate discriminant $\Delta(f)$ is a factor of the squarefreed iterated resultant. In fact, we prove a much stronger result. To be specific, we find a polynomial $\Hp(f,[x_1,\ldots,x_n])$, which is a factor of the squarefreed iterated resultant, and show that the multivariate discriminant $\Delta(f)$ is a factor of $\Hp(f,[x_1,\ldots,x_n])$. Furthermore, we conjecture that
$$\Hp(f,[x_1,\ldots,x_n])=\Delta(f)\,(*)$$
 holds for generic form $f(\xx_n)$, and show that it is true for generic trivariate form (Theorem \ref{thm:main2}).

 It is worthwhile to note that, in practice, iterated resultants appear frequently in CAD (Cylindrical Algebraic Decomposition) \cite{collins1}, and our proof is based on the results about CAD we established in a series works \cite{Han2012Proving,han2014constructing,han2015open}. Moreover, the polynomial $\Hp(f,[x_1,\ldots,x_n])$ is exactly the Han's projection polynomial for producing open weak CAD, and we have already proved that $(*)$ holds for generic quadratic form $f(\xx_n)$.

The structure of this paper is as follows. In Section 2, we give the formal statement of the main results. In Section 3 and Section 4, we give the proof of Theorem \ref{thm:main} and Theorem \ref{thm:main2}, respectively.
\section{Formal statement of the main results}
In this section, we state the main result of this paper. We hereby introduce some notations first.

If not specified, for a positive integer $n$, let $\xx_n$ be the list of variable $(x_1,\dots,x_n)$.
The polynomial rings $\ZZ[\xx_n]$, $\RR[\xx_n]$ appear in the following definitions can be replaced by $\mathcal{R}[\xx_n]$, where $\mathcal{R}$ is a UFD.
\begin{definition}\label{de:sqrfree}
Suppose $h$$\in \ZZ[\xx_n]$ can be factorized in $\ZZ[\xx_n]$ as:
  $$h=a{h_1}^{i_1}{h_2}^{i_2}\ldots {h_m}^{i_m},$$%l_{1}^{2j_1+1}l_2^{2j_2+1}\ldots l_t^{2j_t+1},$$
where $a\in\ZZ$, $h_i (i=1,\ldots,m)$ are pairwise different irreducible primitive polynomials with positive leading coefficient (under a suitable ordering) and positive degree in $\ZZ[\xx_n]$. Define
  $$\sqrfree(h)=\prod_{i=1}^m{h_i}.$$
  If $h$ is a constant, let $\sqrfree(h)=1.$
\end{definition}
\begin{definition}
The {\em level} of $f\in \ZZ[\xx_n]$ is the largest $j$ such that $\deg(f,{x_j})>0$ where $\deg(f,x_j)$ is the degree of $f$ with respect to $x_j$.
Let $f(\xx_n)\in \R[\xx_n]$, say
$f(\xx_n)=\sum_{i=0}^lc_ix_{n}^{i}, c_l \not\equiv 0,$
where $c_i (i=0,\ldots,l)$ is an element of $\R[\xx_{n-1}]$. Then, the {\em leading coefficient} of $f(\xx_n)$ with respect to $x_{n}$ is $c_l$ and denoted by $\lc(f,x_{n}).$
\end{definition}
\begin{definition} \label{def:brown projection}\cite{McCallum1,McCallum2,brown}
        For a given polynomial $F\in \ZZ[\xx_n]$, if $F$ is with level $n$,
        the iterated resultant projection operator $\Bproj$ is defined as
        $$\Bproj(F,[x_{n}])=\Res(\sqrfree(F),\frac{\partial (\sqrfree(F))}{\partial x_{n}}, x_{n}),$$
        where ``\Res" means the Sylvester resultant. Otherwise $\Bproj(F,[x_{n}])=F$.

        Define
        \begin{align*}
                &\Bproj(F,[x_{n},x_{n-1},\ldots, x_i])\\
                =&\Bproj(\Bproj(F,[x_{n},x_{n-1},\ldots,x_{i+1}]),[x_i]).
        \end{align*}
\end{definition}
We can regard $\Bproj(F,[x_{n},x_{n-1},\ldots, x_i])$ $(i=1,\ldots,n)$ as the squarefreed iterated resultants of a given polynomial $F$.
\begin{definition}
Let $F(x) =c_lx^l+\cdots+c_0\in \R[x]$ with $c_l \neq 0.$
The {\em discriminant} of $F(x)$ is
$$\discrim(F,x)=c_l^{2l-2}\prod_{i<j}{(z_i-z_j)}^2,$$
where $z_i$ $(i=1,\ldots,l)$ are the complex roots of the equation $F(x)=0$.

We have the following well-known relationship,
\begin{equation}\label{eqn:resdis}
c_l\discrim(F,x)=(-1)^{\frac{l(l-1)}{2}}\Res(F,\frac{\partial}{\partial x}F,x).
\end{equation}
\end{definition}

\begin{definition}\cite{han2014constructing}\label{def:hp}[Han's projection operator for open weak CAD]
Let $F\in \ZZ[x_1,\ldots,x_n]$, for a given $m (1\le m\le n),$ denote $[\yy]=[y_1,\dots,y_{m}]$ where $ y_i \in \{x_1,\dots,x_n\}$ $(1\le i\le m)$ and $y_i\neq y_j$ ($i\neq j$). For $1\le i\le m$,
$\Hproj(F,[\yy],y_i)$ and $\Hproj(F,[\yy])$ are defined recursively as follows.
        \begin{align*}
                 \Bproj(F,[y_i])     &=\Res(\sqrfree(F),\frac{\partial (\sqrfree(F))}{\partial y_{i}}, y_{i}),\\
                \Hproj(F,[\yy],y_i) &=\Bproj(\Hproj(F,[\hat{\yy]}_i),[y_i]),\\
                \Hproj(F,[\yy])     &=\gcd(\Hproj(F,[\yy],y_1),\ldots,\Hproj(F,[\yy],y_m)),\\
                \Hproj(F,[~])       &=F,
    \end{align*}
where $\hat{[\yy]}_i=[y_1,\ldots,y_{i-1},y_{i+1},\ldots,y_{m}]$.

It is clear that $\Hproj(F,[x_{n},x_{n-1},\ldots, x_i])$ is a factor of $\Bproj(F,[x_{n},x_{n-1},\ldots, x_i])$.
\end{definition}
\begin{example}
 For any polynomial $F\in \ZZ[x_1,\ldots,x_n]$, we have \[\Hproj(F,\left[x_{1},x_{2}\right]  )=\gcd\left(\Bproj(\Bproj(F,[x_{2}%
]),[x_{1}]),\Bproj(\Bproj(F,[x_{1}]),[x_{2}])\right).
\]
For the generic trivariate quadratic form $F(x,y,z)=ax^2+bxy+cy^2+dxz+eyz+fz^2$, we have
$$\Hproj(F,[x,y,z])=4acf-ae^2-b^2f+bde-cd^2.$$
\end{example}

\begin{definition}
Let $\bm{\alpha}=(\alpha_1,\ldots,\alpha_n)$, $|\bm{\alpha}|=\sum_{i=1}^n \alpha_i$, $\bm {x}^{\bm{\alpha}}=\prod_{i=1}^n x_i^{\alpha_i}$, $\{\bm{C}_{\bm{\alpha}}\}=\{C_{\bm{\alpha}}||\bm{\alpha}|=d\}$, and $N=(\begin{subarray}{c}n+d-1\\n-1 \end{subarray})$. A ``generic'' form $f$ in $n$ variables with degree $d$ is defined as follows,
$$f(\xx_n,\bm{C}_{\bm{\alpha}})=\sum_{|\bm{\alpha}|=d}C_{\bm{\alpha}}\bm {x}^{\bm{\alpha}},$$
where we regard the coefficients of form $f$ as parameters. In some references, $f$ is called a ``universal'' form. For convenience, we will simply write $f(\xx_n,\bm{C}_{\bm{\alpha}})$ as $f(\xx_n)$.
\end{definition}

\begin{definition}
The multivariate discriminant $\Delta(f)$ of the generic form $f(\xx_n)$ is an irreducible polynomial of degree $n(d-1)^{n-1}$ in $\ZZ[\bm{C}_{\bm{\alpha}}]$ satisfying
$$\Delta(f)=0\Longleftrightarrow \exists \bm{u}_n\in\C^n\backslash\{0\}, \nabla_{\xx_n}f(\bm{u}_n)=0.$$
In fact, $\Delta(f)$ is even irreducible in $\C[\bm{C}_{\bm{\alpha}}]$. We refer the reader to \cite{israel1994gel} for more details about multivariate discriminant.
\end{definition}

The main results of this paper are the following theorems.
\begin{theorem}\label{thm:main}
For generic form $f(\xx_n)$ with even degree $d$, $\Delta(f)$ is an irreducible factor of $\Hproj(f,[x_n,\ldots,x_1])$ in $\ZZ[\bm{C}_{\bm{\alpha}}]$. In particular, $\Delta(f)$ is an irreducible factor of the squarefreed iterated resultant $\Bproj(f,[x_n,\ldots,x_1])$ in $\ZZ[\bm{C}_{\bm{\alpha}}]$.
\end{theorem}
\begin{remark}
  Theorem \ref{thm:main} does not hold for all forms $F$ in $n$ variables with even degree $d$, for example, let $n=3$, $d=2$, and
  $$F=xy+y^2+xz+yz+kz^2,$$
  then $\Hproj(F,[x,y,z])=1,$ and $\Delta(F)=-k$. In this case, $\Hproj(F,[x,y,z])|\Delta(F)$.

  We also note that, by definition, while $\Delta(F)$ may vanish identically, $\Hproj(f,[x_n,\ldots,x_1])$ is always non-zero.
\end{remark}
\begin{theorem}\label{thm:main2}
  For generic form $f(x,y,z)$ in three variables with degree $d$, we have
  $$\Hp(f,[x,y,z])=\Delta_{x,y,z}(f).$$
\end{theorem}

In \cite{han2015open}, we show that $\Hproj(f,[x_n,\ldots,x_1])=\Delta(f)$ is true for generic quadratic form $f$ in $n$ variables. It is reasonable to guess that the following conjecture is true.
  \begin{conjecture}For generic form $f(x_1,\ldots,x_n)$, we have
  $$\Hproj(f,[x_n,\ldots,x_1])=\Delta(f).$$
  \end{conjecture}

\section{The proof of Theorem \ref{thm:main}}
In order to prove Theorem \ref{thm:main}, we first introduce some notations and results.

\begin{definition}
For an even integer $d$, let $P_{n,d}$ be the cone of nonnegative forms in $\RR^n$ of degree $d$, that is
$$P_{n,d}=\{C_{N}\in\RR^N|\forall X_n\in\RR^n,f(X_n,C_{N})\ge0\},$$
recall that $N=(\begin{subarray}{c}n+d-1\\n-1 \end{subarray})$. %\Card(\{\bm{C}_{\bm{\alpha}}\})$.
\end{definition}
\begin{remark}\label{rem:boundary}
It is well known that $P_{n,d}$ is a nonempty closed convex cone in $\RR^N$ of full dimension. If $C_{N}\in \partial P_{n,d}$, then $\forall X_n\in\RR^n,f(X_n,C_{N})\ge0$, and for any open neighborhood $S\subseteq \RR^N$ of $C_N$, $f(\xx_n,\bm{C}_{\bm{\alpha}})$ is not positive semi-definite on $\RR^n\times S$.
\end{remark}
\begin{definition}
For a polynomial $g\in \ZZ[\xx_n]$, we define
$$V_{\CC}(g)=\{X_n|X_n\in \CC^n, g(X_n)=0\}, V_{\RR}(g)=\{X_n|X_n\in \RR^n, g(X_n)=0\}.$$
It is clear that $V_{\RR}(g)\subseteq V_{\CC}(g)$.
\end{definition}

The following lemma establishes a connection between $\partial P_{n,d}$ and $V_{\CC}(\Delta(f))$.
\begin{lemma}\cite{nie2012discriminants} \label{thm:nie2012}
For an even integer $d$, the Zariski closure of the boundary $\partial P_{n,d}$ is $V_{\CC}(\Delta(f))$, which is an irreducible hypersurface of degree $n(d-1)^{n-1}$.
\end{lemma}

\begin{definition}\label{def:weakopendeli} \cite{han2015open}
Let $S$ be an open set of $\RR^j$. The polynomial $F(\xx_n)$ is said to be {\em open weak delineable} on $S$ if, for any maximal open connected set $U\subseteq \RR^n$ defined by $F(\xx_n)\neq0$, we have
$(S\times \RR^{n-j})\bigcap U\ne \emptyset \Longleftrightarrow \forall\va\in S, (\va\times \RR^{n-j})\bigcap U\ne \emptyset.$
\end{definition}

The following lemma shows an important geometry property of the projection operator $\Hp$.
\begin{lemma} \label{thm:open delineable} \cite{han2015open}
For generic form $f(\xx_n)$, there exists a nonzero polynomial $h\in\ZZ[\bm{C}_{\bm{\alpha}}]$, such that for any open connected set $S$ of
$$\RR^{N}\backslash V_{\RR}(\Hproj(f,[x_n,\dots,x_1]),$$
$f$ is open weak delineable on $S\backslash V_{\RR}(h)$.
\end{lemma}
%As a direct corollary, we have
The following lemma is a well known result, we state it without proof.
\begin{lemma}\label{lem:div}
  Let $F,G\in\ZZ[x_1,\ldots,x_n]$, if $F$ is a factor of $G$ in $\CC[x_1,\ldots,x_n]$, then $F$ is a factor of $G$ in $\QQ[x_1,\ldots,x_n]$. In particular, when $F$ is a primitive polynomial, $F$ is a factor of $G$ in $\ZZ[x_1,\ldots,x_n]$.
\end{lemma}

Now we can prove Theorem \ref{thm:main}.

\begin{prof}
Let $S$ be any open connected set of $\RR^{N}\backslash V_{\RR}(\Hproj(f,[x_n,\dots,x_1]))$. According to Lemma \ref{thm:open delineable}, $f$ is open weak delineable on the open connected set $S'=S\backslash V_{\RR}(h)$ for some nonzero polynomial $h\in\ZZ[\bm{C}_{\bm{\alpha}}]$.

We claim that for any $\vb_0\in S'$,
$$\forall X_n\in\RR^n, f(X_n,\vb_0)\ge0\Longleftrightarrow \forall \vb\in S',\forall X_n\in\RR^n, f(X_n,\vb)\ge0.$$
 Otherwise, there exists $\vb_1\in S'$, $X_n'\in \RR^n$, such that $f(X_n',\vb_1)<0$. Let $U\subseteq \RR^{N+n}$ be the maximal open connected set defined by $f(\xx_n)\neq0$ containing point $(X_n',\vb_1)$, it is clear that $f(U)<0$. By the definition of open weak delineable,
 $$(\vb_1\times \RR^{N})\bigcap U\ne \emptyset\Longrightarrow(S'\times \RR^{N})\bigcap U\ne \emptyset\Longrightarrow (\vb_0\times \RR^{N})\bigcap U\ne \emptyset,$$
 which means that $f(\xx_n,\vb_0)$ is not positive semi-definite, a contradiction.

We next prove that if $f$ is not positive semi-definite on any point in $S'$, then $f$ is not positive semi-definite on any point $\vb\in S$. Otherwise, we assume that $f(\xx_n,\vb)$ is positive semi-definite. In this case, $S\cap int(P_{n,d})$ is a non-empty open set, since $P_{n,d}$ is a convex cone of full dimension. As a result, $S'\cap int(P_{n,d})=(S\backslash V_{\RR}(h))\cap int(P_{n,d})$ is a non-empty set, and $f$ is positive semi-definite on it, a contradiction.

It is clear that if $f$ is positive semi-definite on $S'$, $f$ is positive semi-definite on $S$, since $P_{n,d}$ is closed.

Therefore, for any $\vb_0\in S$,
$$\forall X_n\in\RR^n, f(X_n,\vb_0)\ge0\Longleftrightarrow \forall \vb\in S,\forall X_n\in\RR^n, f(X_n,\vb)\ge0.$$
By Remark \ref{rem:boundary}, $S \bigcap \partial P_{n,d}=\emptyset.$ Thus, %for any $1\le t\le n$,
$$\partial P_{n,d}\subseteq V_{\RR} (\Hproj(f,[x_n,\dots,x_1])),$$
and their corresponding Zariski closures
$$\Zar(\partial P_{n,d})\subseteq \Zar(V_{\RR} (\Hproj(f,[x_n,\ldots,x_1])))\subseteq V_{\CC}(\Hproj(f,[x_n,\ldots,x_1])).$$
By Lemma \ref{thm:nie2012}, $\Zar(\partial P_{n,d})=V_{\CC}(\Delta(f))$.
Since $\Delta(f)$ is irreducible in $\CC[\bm{C}_{\bm{\alpha}}]$, by Hilbert's Nullstellensatz, $$\Delta(f)|\Hproj(f,[x_n,\ldots,x_1]) \mbox{ in } \CC[\bm{C}_{\bm{\alpha}}].$$
Notice that $\Delta(f)$ and $\Hproj(f,[x_n,\ldots,x_1])$ are polynomials in $\ZZ[\bm{C}_{\bm{\alpha}}]$ and $\Delta(f)$ is irreducible, by Lemma \ref{lem:div}, $\Delta(f)$ is an irreducible factor of $\Hproj(f,[x_n,\ldots,x_1])$ in $\ZZ[\bm{C}_{\bm{\alpha}}]$.
\end{prof}

\section{The proof of Theorem \ref{thm:main2}}
In order to prove Theorem \ref{thm:main2}, we first introduce some notations and results.

\begin{definition}(Multipolynomial Resultants)
  For generic forms
  $$F_i=\sum_{|\bm{\alpha}_i|=d_i}\bm{C}_{i,\bm{\alpha}_i}\xx^{\bm{\alpha}_i},\, (i=1,\ldots,n)$$
  the resultant of $F_i$ is the unique polynomial $\Res\in\ZZ[\bm{C}_{i,\alpha_i}]$ (up to a constant) which has the following properties:

  (1) $F_i=0$ has a nontrivial solution over $\CC$ $\Longleftrightarrow$ $\Res(F_1,\ldots,F_n)=0$.

  (2) $\Res(x_1^{d_1},\ldots,x_n^{d_n})=1$.

  (3) $\Res$ is irreducible in $\CC[\bm{C}_{i,\bm{\alpha}_i}]$.

\end{definition}
We refer the reader to \cite{israel1994gel,cox2006using} for more details about multipolynomial resultants.

\begin{definition}
For a polynomial $F(x,y,z)\in \ZZ[x,y,z]$ and a given integer $i$, $\delta_{z,z'}^i(F)$ is defined as
$$\delta_{z,z'}^i(F)=\sum_{k\ge i}\frac{1}{k!}(z'-z)^{k-i}\frac{\partial^k F(x,y,z)}{\partial z^k}.$$
It is easy to verify that $\delta_{z,z'}^i(F)$ satisfies
$$F(x,y,z')=\sum_{j=0}^{i-1}(z'-z)^j\frac{\partial^j F(x,y,z)}{\partial z^j}+(z'-z)^i\delta_{z,z'}^i(F).$$

\end{definition}

The following lemma is a restatement of the main result in \cite{buse2009explicit}.
\begin{lemma}\cite{buse2009explicit}\label{lem:buse}%Assuming that $d\ge4$,
  For generic trivariate form $f(x,y,z)$ with degree $d\ge3$, we have the decomposition in irreducible factors in $\ZZ[\bm{C}_{\bm{\alpha}}]$ (up to a constant),
    $$\discrim(\discrim(f(1,y,z),z),y)=C_{0,0,d}\Delta_{x,y,z}(f)a_{z,y}(f)^3b_{z,y}(f)^2,$$
    where $a_{y,z}$ and $b_{y,z}$ are two irreducible polynomials in $\ZZ[\bm{C}_{\bm{\alpha}}]$, such that
  $$\Res_{x,y,z}(f,\frac{\partial f}{\partial z}, \frac{\partial^2 f}{\partial z^2})=2^{d(d-1)}C_{0,0,d}^2a_{z,y}(f),$$
  and
  $$\Res_{x,y,z,z'}(f,\frac{\partial f}{\partial z}, \delta_{z,z'}^2(f),\delta_{z,z'}^2(\frac{\partial f}{\partial z})-2\delta_{z,z'}^3(f))=C_{0,0,d}^{2d(d-1)-6}b_{z,y}(f)^2,$$
  for $d\ge4$, $b_{z,y}(f)=1$ for $d=3$.
    \end{lemma}
\begin{remark}
  The main theorem in \cite{buse2009explicit} assumes $d\ge4$, but its proof is also valid for the case $d=3$.
\end{remark}
%  2) if $d=3$,   $$\discrim(\discrim(f(1,y,z),z),y)=C_{0,0,d}\Delta_{x,y,z}(f)a_{z,y}(f)^3,$$
Now, we can prove Theorem \ref{thm:main2}.

\begin{prof}
It is easy to check by hand (or with a computer) for the cases $d = 1,2$.

For $d\ge3$, it suffices to prove that
$$\sqrfree(\Hp(f,[y,z]))=\Delta(f)\Delta(f(0,y,z))x,$$
since $$\Hp(f,[y,z],x)=\Bp(\Hp(f,[y,z]),x),$$ $$\gcd(\Delta_{y,z}(f(0,y,z)),\Delta_{x,z}(f(x,0,z)),\Delta_{x,y}(f(x,y,0)))=1,$$
and
$$\Hp(f,[x,y,z])=\gcd(\Hp(f,[y,z],x),\Hp(f,[x,z],y),\Hp(f,[y,x],z)).$$

By definition,
  $$\Hp(f,[y,z])=\gcd(\Hp(f,[y],z),\Hp(f,[z],y))=\gcd(\Bp(f,[y,z]),\Bp(f,[z,y])).$$
 According to Lemma \ref{lem:buse}, % and Equation (\ref{eqn:resdis}) For simplicity, for $d=3$, we denote $b_{z,y}(f)=a_{z,y}(f)=1$,
  \begin{align*}
\sqrfree(\Bp(f,[z,y]))=&x\sqrfree(\Bp(f(1,y,z),[z,y]))\\
=&x\sqrfree(\Bp(C_{0,0,d}\discrim(f(1,y,z),z),[y]))\\
=&x\sqrfree(C_{0,0,d}\discrim(\discrim(f(1,y,z),z),y)\lc(\discrim(f(1,y,z),z),y))\\
=&x\sqrfree(C_{0,0,d}\Delta_{x,y,z}(f)a_{z,y}(f)b_{z,y}(f)\discrim(f(0,1,z),z))\\
=&x\sqrfree(C_{0,0,d}\Delta_{x,y,z}(f)a_{z,y}(f)b_{z,y}(f)\Delta_{y,z}(f(0,y,z))),
\end{align*}
where we use the fact that
$$\Delta_{y,z}(f(0,y,z))=\discrim(f(0,1,z),z).$$
Similarly,
$$\sqrfree(\Bp(f,[y,z]))=x\sqrfree(C_{0,d,0}\Delta_{x,y,z}(f)a_{y,z}(f)b_{y,z}(f)\Delta_{y,z}(f(0,y,z))).$$

When $d=3$, $b_{y,z}(f)=b_{z,y}(f)=1$, and it is easy to check with a computer that $a_{y,z}(f),a_{z,y}(f)$ are two different irreducible polynomials. Thus, $$\sqrfree(\Hp(f,[y,z]))=\gcd(\Bp(f,[y,z]),\Bp(f,[z,y]))=\Delta(f)\Delta(f(0,y,z))x.$$%$b_{y,z}(f),b_{z,y}(f)$ are two constants,

When $d\ge4$, let $L_1(f)=\{a_{z,y}(f)$, $b_{z,y}(f)\}$, $L_2(f)=\{a_{y,z}(f)$, $b_{y,z}(f)\}$ be two polynomial sets.
According to Lemma \ref{lem:buse}, the polynomials in $L_1(f)\bigcup L_2(f)$ are irreducible. In order to prove $\sqrfree(\Hp(f,[y,z]))=\Delta(f)\Delta(f(0,y,z))x$, we only need to show that $L_1(f)\bigcap L_2(f)=\emptyset$, and it suffices to find a polynomial $F(x,y,z,w)\in\ZZ[x,y,z,w]$, which is homogenous in $x,y,z$, such that $L_1(F)\bigcap L_2(F)=\emptyset$. % $a_{z,y}(F)\notin L_2(F)$.%We first prove that $a_{z,y}(f)\notin L_2(f)$.
%%$$\sqrfree(\Hp(f,[y,z]))=\Delta_{x,y,z}(f)\Delta_{y,z}(f(0,y,z))x,$$

Let $F(x,y,z,w)=z^d+wzx^{d-1}+y^d$. We first prove that $a_{y,z}(F)=b_{y,z}(F)=0$.
It is easy to see that $\Res_{x,y,z}(F,\frac{\partial F}{\partial y}, \frac{\partial^2 F}{\partial y^2})=a_{y,z}(F)=0$, since
 $\{F=0,\frac{\partial F}{\partial y}=0, \frac{\partial^2 F}{\partial y^2}=0\}=\{z^d+wzx^{d-1}+y^d=0,dy^{d-1}=0,d(d-1)y^{d-1}=0\}$ always has a nontrivial solution $(1,0,0)$.

 By definition of $\delta_{y,y'}^iF$,
$$\delta_{y,y'}^2(\frac{\partial F}{\partial y})-2\delta_{y,y'}^3(F)|_{y=y'=0}=0,$$
$$\delta_{y,y'}^2(F)|_{y=y'=0}=0,$$
thus,
$\{F=0,\frac{\partial F}{\partial y}=0, \delta_{y,y'}^2(F)=0,\delta_{y,y'}^2(\frac{\partial F}{\partial y})-2\delta_{y,y'}^3(F)=0\}$
always has a nontrivial solution $(x,y,z,y')=(1,0,0,0)$, and
$$\Res_{x,y,z,y'}(F,\frac{\partial F}{\partial y}, \delta_{y,y'}^2(F),\delta_{y,y'}^2(\frac{\partial F}{\partial y})-2\delta_{y,y'}^3(F))=0, b_{y,z}(F)=0.$$

Now, we prove that $w|a_{z,y}(F)\neq0,$ and $w|b_{z,y}(F)\neq0$.

 It is clear that $\{F=0,\frac{\partial F}{\partial z}=0, \frac{\partial^2 F}{\partial z^2}=0\}=\{z^d+wzx^{d-1}+y^d=0,dz^{d-1}+wx^{d-1}=0,d(d-1)z^{d-2}=0\}$ has a nontrivial solution if and only if $w=0$. Thus, $w|a_{z,y}(F)\neq0.$

  By definition of $\delta_{z,z'}^iF$,
  \begin{align*}
  \delta_{z,z'}^2(F)&=\frac{z'^d+wz'x^{d-1}+y^d-z^d-wzx^{d-1}-y^d-(z'-z)(dz^{d-1}+wx^{d-1})}{(z'-z)^2}\\
  &=\frac{z'^d-z^d-d(z'-z)z^{d-1}}{(z'-z)^2}\\
  &=\sum_{k\ge 2}\frac{1}{k!}(z'-z)^{k-i}\frac{\partial^k F(x,y,z)}{\partial z^k},
  \end{align*}
\begin{align*}
\delta_{z,z'}^2(\frac{\partial F}{\partial z})-2\delta_{z,z'}^3(F)=&\frac{dz'^{d-1}-dz^{d-1}-d(d-1)(z'-z)z^{d-2}}{(z'-z)^2}-\\
&2\cdot\frac{z'^d-z^d-d(z'-z)z^{d-1}-d(d-1)(z'-z)^2z^{d-2}}{(z'-z)^3}.
\end{align*}
  Suppose that
  $$L=\{\delta_{z,z'}^2(F)=0,\delta_{z,z'}^2(\frac{\partial F}{\partial z})-2\delta_{z,z'}^3(F)=0\}$$
 has a nontrivial solution $(z,z')=(p,p+q)$.
 We must have $p\neq0$, $q\neq 0$, since
 $$\delta_{z,z'}^2(F)|_{z=0,z'=q}=q^{d-2},\, \delta_{z,z'}^2(F)|_{z=z'=p}=\frac{d(d-1)p^{d-2}}{2}.$$

 Notice that the polynomials in $L$ are homogenous, we may assume $p=1$. Since
 \begin{align*}
 &(dz'^{d-1}-dz^{d-1}-d(d-1)(z'-z)z^{d-2})(z'-z)-\\
 &2(z'^d-z^d-d(z'-z)z^{d-1}-d(d-1)(z'-z)^2z^{d-2})|_{z=1,z'=1+q}\\
 =&(d(1+q)^{d-1}-d-d(d-1)q)q-2((1+q)^d-1-dq-d(d-1)q^2)\\
 =&d(d-1)q^2+dq+2+d(1+q)^{d-1}q-2(1+q)^d,
 \end{align*}
 $q$ is a nonzero solution of $\{(q+1)^d=1+dq,2(1+q)^d=d(d-1)q^2+dq+2+d(1+q)^{d-1}q\}$. From the first equation, we know that $q\neq-1$. Substitute the equality $(q+1)^d=1+dq$ into the second equation, we have
 $$2(1+q)^d=d(d-1)q^2+dq+2+d(1+q)^{d-1}q,$$
 $$\Longleftrightarrow2(1+dq)=d(d-1)q^2+dq+2+\frac{dq(1+dq)}{1+q},$$
 $$\Longleftrightarrow dq(1+q)=d(d-1)q^2(1+q)+dq(1+dq),$$
 $$\Longleftrightarrow 1+q=(d-1)q(1+q)+1+dq,$$
 $$\Longleftrightarrow 1=(d-1)(1+q)+d\Longleftrightarrow q=-2,$$
 and
 $$(q+1)^d=1+dq\Longrightarrow (-1)^d=1-2d \Longrightarrow d=1,$$
 a contradiction.

 Hence $L$ only has a nontrivial solution $(z,z')=(0,0)$, and
 $$\{F|_{z=0}=0,\frac{\partial F}{\partial z}|_{z=0}=0\}=\{y^d=0,wx^{d-1}=0\}$$
  only has a nontrivial solution if and only if $w=0$, which implies that
  $$L\bigcup \{F=0,\frac{\partial F}{\partial z}=0\}$$
  has a nontrivial solution if and only if $w=0$, and $w|b_{z,y}(F)\neq0$.

Therefore $\{a_{z,y}(F),b_{z,y}(F)\}\bigcap\{a_{y,z}(F),b_{y,z}(F)\}=\emptyset$, and we are done.
\end{prof}

\textbf{Acknowledgements}{\rm \quad The author would like to thank his advisor Gang Tian for constant encouragement and several useful comments on an earlier version of this paper. The author would want to convey his gratitude to Laurent Bus{\'e} for his help on illustrating some details in their paper. The author would also want to thank China Scholarship Council (CSC) and ``Training, Research and Motion" (TRAM) network for supporting him visiting Princeton University. We thank the referees for their time and comments.}
%\begin{thebibliography}{99}

\end{document}